\documentclass[12pt,leqno]{article}
\usepackage{amssymb,amsfonts,amsmath,amsthm,amscd,mathrsfs}
\usepackage{PSFigure,psfrag}
\usepackage{graphics}
\usepackage{graphicx}
\def\IP{{\mathbb P}}

\def\IN{{\mathbb N}}

\def\IZ{{\mathbb Z}}

\def\n{\noindent}
\def\dsl{\textstyle\sum\limits}

\def\dis{\displaystyle}
\def\o{\omega}

\def\ov{\overline}
\def\ve{\varepsilon}
\def\f{\footnotesize}
\def\r{\rightarrow}

\def\wh{\widehat}
\def\wt{\widetilde}

\def\cA{{\cal A}}

\def\cT{{\cal T}}

\def\cI{{\cal I}}

\def\cH{{\cal H}}

\def\cV{{\cal V}}
\def\cW{{\cal W}}

\dimendef\dimen=0

\newtheorem{theorem}{Theorem}[section]
\newtheorem{lemma}[theorem]{Lemma}

\newtheorem{proposition}[theorem]{Proposition}
\newtheorem{remark}[theorem]{Remark}

\begin{document}

\baselineskip14pt
\noindent

\begin{center}
{\bf CONNECTIVITY BOUNDS FOR THE VACANT SET \\ OF RANDOM INTERLACEMENTS}
\end{center}

\vspace{1cm}
\begin{center}
Vladas Sidoravicius\footnote[1]{CWI, Kruislaan 413, NL-1098 SJ, Amsterdam, and IMPA, Estrada Dona Castorina 110, Jardim 

 ~~Botanico, CEP 22460-320, Rio de Janeiro, RJ, Brasil.} and Alain-Sol Sznitman\footnote[2]{Departement Mathematik, ETH Z\"urich, CH-8092 Z\"urich, Switzerland.
 

}
\end{center}


\bigskip
\begin{abstract}
The model of random interlacements on $\IZ^d$, $d \ge 3$, was recently introduced in \cite{Szni07a}. A non-negative parameter $u$ parametrizes the density of random interlacements on $\IZ^d$. In the present note we investigate connectivity properties of the vacant set left by random interlacements at level $u$, in the non-percolative regime $u > u_*$, with $u_*$ the non-degenerate critical parameter for the percolation of the vacant set, see \cite{Szni07a}, \cite{SidoSzni08}. We prove a stretched exponential decay of the connectivity function for the vacant set at level $u$, when $u > u_{**}$, where $u_{**}$ is another critical parameter introduced in \cite{Szni08b}. It is presently an open problem whether $u_{**}$ actually coincides with $u_*$.
\end{abstract}





\newpage

\thispagestyle{empty}
~

\newpage
\setcounter{page}{1}

\setcounter{section}{-1}
\section{Introduction}

In this note we derive stretched exponential bounds on the connectivity function for the vacant set of random interlacements on $\IZ^d$, $d \ge 3$. The model of random interlacements has been introduced in \cite{Szni07a}. It heuristically describes the microscopic structure left in the bulk by random walk on a cylinder with base a large $(d-1)$-dimensional discrete torus, or by random walk on a large $d$-dimensional torus, when the walk respectively runs for times proportional to the square of the number of sites in the base, or to the number of sites in the torus, cf.~\cite{Szni07b}, \cite{Wind08a}. Further extensions to more general graphs can be found in \cite{Teix08b}, \cite{Wind08b}.

\medskip
The bounds presented here pertain to a very specific region of the non-percolative regime of the vacant set of random interlacements. Knowing whether this region actually coincides with the whole non-percolative regime outside the critical point is an important question with direct implications for the asymptotic behavior of the disconnection time of discrete cylinders with bases which become large, cf.~\cite{Szni08b}. The results in the present note do not answer this question, but they indicate that in the case where the two regions differ, there is a marked transition in the decay properties of the connectivity function (roughly from a polynomial to a stretched exponential decay), as one moves across another (and then distinct) critical point, which has been introduced in \cite{Szni08b}.

\medskip
We will now describe the model and state our main result. We refer to Section 1 for precise definitions. Random interlacements consists of a cloud of paths which constitute a Poisson point process on the space of doubly infinite $\IZ^d$-valued trajectories modulo time-shift, tending to infinity at positive and negative infinite times. A non-negative parameter $u$ plays the role of a multiplicative factor of the intensity measure of this Poisson point process. In a standard fashion one constructs on the same probability space $(\Omega, \cA, \IP)$, see below (\ref{1.6}), the whole family $\cI^u$, $u \ge 0$, of random interlacements at level $u \ge 0$, cf.~(\ref{1.10}). They come as traces on $\IZ^d$ of the cloud of trajectories modulo time-shift with labels at most $u$. The subsets $\cI^u$ increase with $u$  and for $u > 0$ are random connected subsets of $\IZ^d$, ergodic under space translations, cf.~Theorem 2.1, Corollary 2,3 of \cite{Szni07a}. The complement $\cV^u$ of $\cI^u$ in $\IZ^d$ is the so-called vacant set at level $u$. It is known from Theorem 3.5 of \cite{Szni07a} and Theorem 3.4 of \cite{SidoSzni08} that there is a non-degenerate critical value $u_* \in (0,\infty)$ such that 
\begin{equation}\label{0.1}
\begin{array}{rl}
{\rm i)} &\mbox{for $u > u_*$, $\IP$-a.s. all connected components of $\cV^u$ are finite},
\\[1ex]
{\rm ii)} &\mbox{for $u < u_*$, $\IP$-a.s. there exists an infinite component in $\cV^u$}.
\end{array}
\end{equation}

\n
It is also known from \cite{Teix08a}, that when $u$ is such that $\cV^u$ percolates (i.e. there is $\IP$-a.s. an infinite connected component in $\cV^u$), the infinite cluster is almost surely unique. It is presently unknown whether $\cV^{u_*}$ percolates or not. Another critical point $u_{**} \in [u_*,\infty)$ has been introduced in \cite{Szni08b}:
\begin{equation}\label{0.2}
\begin{array}{l}
u_{**}  = \inf\{ u  \ge 0; \alpha (u) > 0\}, \;\mbox{where}
\\[1ex]
\alpha(u) = \sup\{ \alpha \ge 0; \lim\limits_{L \r \infty} L^\alpha\, \IP[B(0,L) \stackrel{\cV^u}{\longleftrightarrow} S(0,2L)] = 0\}, \;\mbox{for $u \ge 0$}, 
\end{array}
\end{equation}

\n
here the event under the probability refers to the existence of a nearest neighbor path in $\cV^u$ joining $B(0,L)$, the closed ball or radius $L$ and center $0$ for the $\ell^\infty$-distance, with $S(0,2L)$, the sphere with radius $2L$ and center $0$ for the same distance. The supremum in the second line of (\ref{0.2}) is by convention equal to zero when the set is empty, (this is for instance the case when $u < u_*$). The critical parameter $u_{**}$ explicitly enters the upper bound on the disconnection time of discrete cylinders derived in Corollary 4.6 of \cite{Szni08b}. It is an important question whether in fact $u_*$ and $u_{**}$ coincide. The present work does not address this question but shows that the probability, which appears in the second line of (\ref{0.2}), has a stretched exponential decay in $L$ for $u > u_{**}$. Our main result is

\begin{theorem}\label{theo0.1}  ~

\medskip
For $u > u_{**}$, the connectivity function in the vacant set at level $u$ has stretched exponential decay. Namely there exist positive constants $\alpha_1,\alpha_2$, and $0< \rho<1$, solely depending on $d$ and $u$, such that with similar notation as in {\rm (\ref{0.2})}:
\begin{equation}\label{0.3}
\IP[0 \stackrel{\cV^u}{\longleftrightarrow} x] \le \alpha_1 \, \exp\{-\alpha_2 |x|^\rho\}, \;\mbox{for all $x \in \IZ^d$}\,.
\end{equation}
\end{theorem}

In case $u_*$ and $u_{**}$ differ, the above theorem forces a sharp transition in the decay properties of the connectivity function of the vacant set as $u$ crosses the value $u_{**}$, cf.~Remark \ref{rem3.1} 1). We also discuss in Remark \ref{rem3.1} 2) a variant of (\ref{0.3}), cf.~(\ref{3.12}), where $\cI^u$ is replaced by its $R$-neighborhood, and instead of assuming $u > u_{**}$, we pick $R$ large enough. Let us mention that when $d=3$, the left-hand side of (\ref{0.3}) does not decay exponentially in $|x|$, cf.~Remark \ref{rem3.1} 3).

\medskip
We now give some comments on the proof of Theorem \ref{theo0.1}. The main difficulty stems from the long range dependence of random interlacements. We use an adaptation of the renormalization and sprinkling technique, which appears in Section 3 of \cite{Szni07a}. The main difference resides in the fact that here we separate the combinatorial complexity bounds from the probabilistic estimates on crossings we derive, see (\ref{2.11}), Lemma \ref{lem2.1} and Proposition \ref{prop2.2}. A similar separation is for instance also used in \cite{Teix08c} for the very percolative regime of small $u$. In our context this separation leads to finer probabilistic estimates in the renormalization scheme, which instead of producing polynomial decay in $L$ for quantities such as $\IP[B(0,L) \stackrel{\cV^u}{\longleftrightarrow} S(0,2L)]$, when $u$ is sufficiently large, cf.~Section 3 of \cite{Szni07a}, yield a stretched exponential decay of such quantities, when $u > u_{**}$.

\medskip
We will now describe the organization of this note.

\medskip
In Section 1 we introduce further notation and recall some facts about random interlacements. In Section 2 we develop the renormalization scheme. The main result is Proposition \ref{prop2.2}. The main consequences of this proposition for the next section appears in Proposition \ref{prop2.5}. In Section 3 we complete the proof of Theorem \ref{theo0.1}. We provide further comments and discuss some extensions in Remark \ref{rem3.1}.

\medskip
Finally let us explain the convention we use for constants. Throughout the text $c$ or $c^\prime$ denote positive constants that solely depend on $d$, with values changing from place to place. Numbered constants $c_0,c_1,\dots$ are fixed and refer to the value pertaining to their first appearance in the text. Dependence of constants on additional parameters appears in the notation, so that for instance $c(u)$ denotes a constant depending on $d$ and $u$.

\section{Notation and a brief review of random interlacements}
\setcounter{equation}{0}

In this section we introduce some notation and recall some facts concerning random interlacements. A more detailed review of random interlacements can also be found in Section 1 of \cite{SidoSzni08}.

We let $|\cdot |$ and $|\cdot |_\infty$ respectively denote the Euclidean and the $\ell^\infty$-distance on $\IZ^d$. Throughout we implicitly assume $d \ge 3$. By finite path we mean a sequence $x_0,x_1,\dots,x_N$ in $\IZ^d$ such that $x_i$ and $x_{i+1}$ are neighbors, i.e. $|x_{i+1} - x_i| = 1$, for $0 \le i < N$. We sometimes simply write path when this causes no confusion. The notation $B(x,r)$ and $S(x,r)$ with $x \in \IZ^d$ and $r \ge 0$, for $|\cdot |_\infty$-balls and spheres is explained below (\ref{0.2}). For $A,B$ subsets of $\IZ^d$, we write $A + B$ for the set of $x + y$ with $x$ in $A$ and $y$ in $B$, and $d(A,B) = \inf\{|x-y|_\infty$; $x \in A, y \in B\}$ for the mutual $\ell^\infty$-distance between $A$ and $B$. When $A$ is a singleton $\{x\}$, we simply write $d(x,A)$. The notation $K \subset \subset \IZ^d$ means that $K$ is a finite subset of $\IZ^d$. Given $U \subseteq \IZ^d$, we denote with $|U|$ the cardinality of $U$, with $\partial U = \{x \in U^c ; \exists y \in U$, $|x-y| = 1\}$ the boundary of $U$ and with $\partial_{\rm int} U = \{x \in U; \exists y \in U^c, |x- y| = 1\}$, the interior boundary of $U$.

\medskip
We denote with $W_+$ the space of nearest neighbor $\IZ^d$-valued trajectories defined for non-negative times and tending to infinity. We let $\cW_+$, $(X_n)_{n \ge 0}$, $(\theta_n)_{n \ge 0}$, stand for the canonical $\sigma$-algebra, the canonical process, and the canonical shift on $W_+$. Since we assume $d \ge 3$, simple random walk is transient on $\IZ^d$, and we denote with $P_x$ the restriction of the canonical law of simple random walk starting at $x \in \IZ^d$, to the set $W_+$, which has full measure. We write $E_x$ for the corresponding expectation. We also define $P_\rho = \sum_{x \in \IZ^d} \rho(x) P_x$, when $\rho$ is a measure on $\IZ^d$ and write $E_\rho$ for the corresponding expectation. Given $U \subseteq \IZ^d$, we let $H_U = \inf\{n \ge 0; X_n \in U\}$, $\wt{H}_U = \inf\{n \ge 1; X_n \in U\}$, and $T_U = \inf\{n \ge 0; X_n \notin U\}$, respectively stand for the entrance time in $U$, the hitting time of $U$, and the exit time of $U$.

\medskip
We denote with $g(\cdot,\cdot)$ the Green function of the walk:
\begin{equation}\label{1.1}
g(x,y) = \dsl_{n \ge 0} P_x [X_n = y], \;\;x,y \in \IZ^d\,.
\end{equation}

\n
It is a symmetric function and due to translation invariance $g(x,y) = g(y-x)$, where $g(y) = g(0,y)$. Given $K \subset \subset \IZ^d$, we write $e_K$ for the equilibrium measure of $K$ and ${\rm cap}(K)$ for its total mass, the so-called capacity of $K$:
\begin{equation}\label{1.2}
e_K(x) = P_x[\wt{H}_K = \infty] \,1_K(x), x \in \IZ^d, \;{\rm cap}(K) = \dsl_{x \in K} P_x [\wt{H}_K = \infty]\,.
\end{equation}
The capacity is subadditive, (this fact easily follows from (\ref{1.2})):
\begin{equation}\label{1.3}
{\rm cap}(K \cup K^\prime) \le {\rm cap}(K) + {\rm cap}(K^\prime), \;\mbox{for $K, K^\prime \subset\subset \IZ^d$}\,.
\end{equation}
Further the probability to enter $K$ can be expressed as:
\begin{equation}\label{1.4}
P_x[H_K < \infty] = \dsl_{y \in K} g(x,y) \,e_K(y), \;\mbox{for $x \in \IZ^d$},
\end{equation}
and one has the bounds, cf. (\ref{1.9}) of \cite{Szni07a}
\begin{equation}\label{1.5}
\begin{array}{l}
\dsl_{y \in K} g(x,y) / \sup\limits_{z \in K} \Big(\dsl_{y \in K} g(z,y)\Big) \le 
\\
P_x [H_K < \infty] \le \dsl_{y \in K} g(x,y) / \inf\limits_{z \in K} \Big(\dsl_{y \in K} g(z,y)\Big), \;\mbox{for $x\in \IZ^d$}\,.
\end{array}
\end{equation}
With classical bounds on the Green function, cf.~\cite{Lawl91}, p.~31, it then follows that
\begin{equation}\label{1.6}
c \,L^{d-2} \le {\rm cap} (B(0,L)) \le c^\prime \,L^{d-2}, \;\mbox{for $L \ge 1$}\,.
\end{equation}

\medskip\n
We now turn to random interlacements. They are defined on a probability space $(\Omega, \cA, \IP)$, where a certain canonical Poisson point process can be constructed, and we refer to (1.16), (1.42) of \cite{Szni07a} or (1.9) - (1.12) of \cite{SidoSzni08} for the precise definition of this probability space. In the present note we will only use the fact one can define on $(\Omega, \cA, \IP)$ families of finite Poisson point processes on $(W_+, \cW_+)$, $\mu_{K,u}(dw)$, $u \ge 0$, $K \subset \subset \IZ^d$, and $\mu_{K,u^\prime,u}(dw)$, $0 \le u^\prime < u$, $K \subset \subset \IZ^d$, so that
\begin{align}
&\mbox{$\mu_{K,u^\prime,u}$ and $\mu_{K,u^\prime}$ are independent with respective intensity measures} \label{1.7}
\\
&\mbox{$(u-u^\prime) \,P_{e_K}$ and $u^\prime P_{e_K}$, for any $0 \le u^\prime < u$ and $K \subset \subset \IZ^d$}, \nonumber
\\[2ex]
&\mbox{$\mu_{K,u} = \mu_{K,u^\prime} + \mu_{K,u^\prime,u}$ for any $0 \le u^\prime < u$ and $K \subset \subset \IZ^d$}. \label{1.8}
\end{align}

\n
Moreover the following compatibility relations hold for $K \subset K^\prime \subset \subset \IZ^d$:
\begin{equation}\label{1.9}
\mu_{K,u} = \dsl^m_{i=0} \delta_{\theta_{H_K}(w_i)} 1 \{H_K(w_i)< \infty\}, \;\mbox{if} \;\mu_{K^\prime,u} = \dsl^m_{i=0} \delta_{w_i},
\end{equation}

\n
together with similar compatibility relations with $\mu_{K,u^\prime,u}$ and $\mu_{K^\prime,u^\prime,u}$ in place of $\mu_{K,u}$ and $\mu_{K^\prime,u}$. We refer for instance to (1.13) - (1.15) of \cite{SidoSzni08}, or to (1.18) - (1.21) and Proposition 1.3 of \cite{Szni07a}, for more details.

\medskip
Given $\o \in \Omega$, the interlacement at level $u \ge 0$ is the random subset of $\IZ^d$ defined for $\omega \in \Omega$ via:
\begin{equation}\label{1.10}
\cI^u(\o) = \textstyle{\bigcup\limits_{K \subset\subset \IZ^d}} \;\;\textstyle{\bigcup\limits_{w \in \, {\rm Supp} \,\mu_{K,u}(\o)}} w(\IN)\,,
\end{equation}
where the notation ${\rm Supp} \,\mu_{K,u}(\o)$ refers to the support of the finite point measure $\mu_{K,u}(\o)(dw)$, and $\IN = \{0,1,\dots\}$. The vacant set at level $u$ is then defined as
\begin{equation}\label{1.11}
\cV^u (\o) = \IZ^d \backslash \cI^u(\o), \;\mbox{for} \;\omega \in \Omega, \, u \ge 0\,.
\end{equation}
One finds that, cf.~(1.54) of \cite{Szni07a}:
\begin{equation}\label{1.12}
\cI^u(\o) \cap K = \textstyle{\bigcup\limits_{w \in \,{\rm Supp} \,\mu_{K^\prime,u}(\omega)}} w(\IN) \cap K, \; \mbox{for} \; K \subset K^\prime \subset \subset \IZ^d, u \ge 0, \omega \in \Omega \,.
\end{equation}
It also follows from (\ref{1.7}) that
\begin{equation}\label{1.13}
\IP [\cV^u \supseteq K] = \exp \{ - u \, {\rm cap} (K)\}, \;\mbox{for all} \;K \subset \subset \IZ^d, \,u \ge 0\,.
\end{equation}

\n
This concludes this short review of random interlacements, which will suffice for the purpose of the present note.

\section{The renormalization scheme}
\setcounter{equation}{0}

In this section we develop the renormalization scheme, which will be the main tool in the derivation of Theorem \ref{theo0.1}. It comes as a variation on the method developed in Section 3 of \cite{Szni07a}, and in particular uses the sprinkling technique of \cite{Szni07a} to control the long range interactions present in the model. The main step comes in Proposition \ref{prop2.2} and its consequences for the proof of Theorem \ref{theo0.1} in the next section appear in Proposition \ref{prop2.5}.

\medskip
We introduce a sequence of length scales
\begin{equation}\label{2.1}
\mbox{$ L_n = L_0 \,\ell_0^n$,  for $n \ge 0$,  where $L_0 \ge 1$ and $\ell_0 \ge 100$ is a multiple of $10$}.
\end{equation}
We organize $\IZ^d$ in a hierarchical fashion with $L_0$ corresponding to the finest scale and $L_1 < L_2 < \dots$ to coarser and coarser scales. To this effect we introduce the set of labels at level $n$:
\begin{equation}\label{2.3}
I_n = \{n\} \times \IZ^d, \,n \ge 0\,,
\end{equation}

\n
and to each $m = (n,i) \in I_n$, with $n \ge 0$, we attach the boxes
\begin{equation}\label{2.4}
C_m = (i L_n + [0,L_n)^d) \cap \IZ^d, \;\wt{C}_m = \textstyle{\bigcup\limits_{m^\prime \in I_n, d(C_{m^\prime}, C_m) \le 1}} C_{m^\prime}\,.
\end{equation}
We write $S_m = \partial_{\rm int} \,C_m$ and $\wt{S}_m = \partial_{\rm int} \wt{C}_m$, for $m \in I_n, n \ge 0$. Given $m \in I_n, n \ge 1$, we consider $\cH_1(m)$, $\cH_2(m) \subseteq I_{n-1}$ defined by
\begin{equation}\label{2.5}
\begin{split}
\cH_1(m) & = \{ \ov{m} \in I_{n-1}; \; C_{\ov{m}} \subseteq C_m \;\mbox{and} \;C_{\ov{m}} \cap S_m \not= \emptyset\}
\\[1ex]
\cH_2(m) & = \Big\{ \ov{m} \in I_{n-1}; \; C_{\ov{m}} \cap  \Big\{z \in \IZ^d; d(z,C_m) = \mbox{\f $\dis\frac{L_{n}}{2}$}\Big\} \not= \emptyset\Big\}\,.
\end{split}
\end{equation}
Note that for all $n \ge 1$, $m \in I_n$:
\begin{equation}\label{2.6}
\begin{array}{l}
\mbox{$\ov{m}_1 \in \cH_1(m)$ and $\ov{m}_2 \in \cH_2(m)$ implies that $\wt{C}_{\ov{m}_1} \cap \wt{C}_{\ov{m}_2} = \emptyset$ and}
\\
\wt{C}_{\ov{m}_1} \cup \wt{C}_{\ov{m}_2} \subseteq \wt{C}_m\,.
\end{array}
\end{equation}

\psfragscanon
\begin{center}
\includegraphics{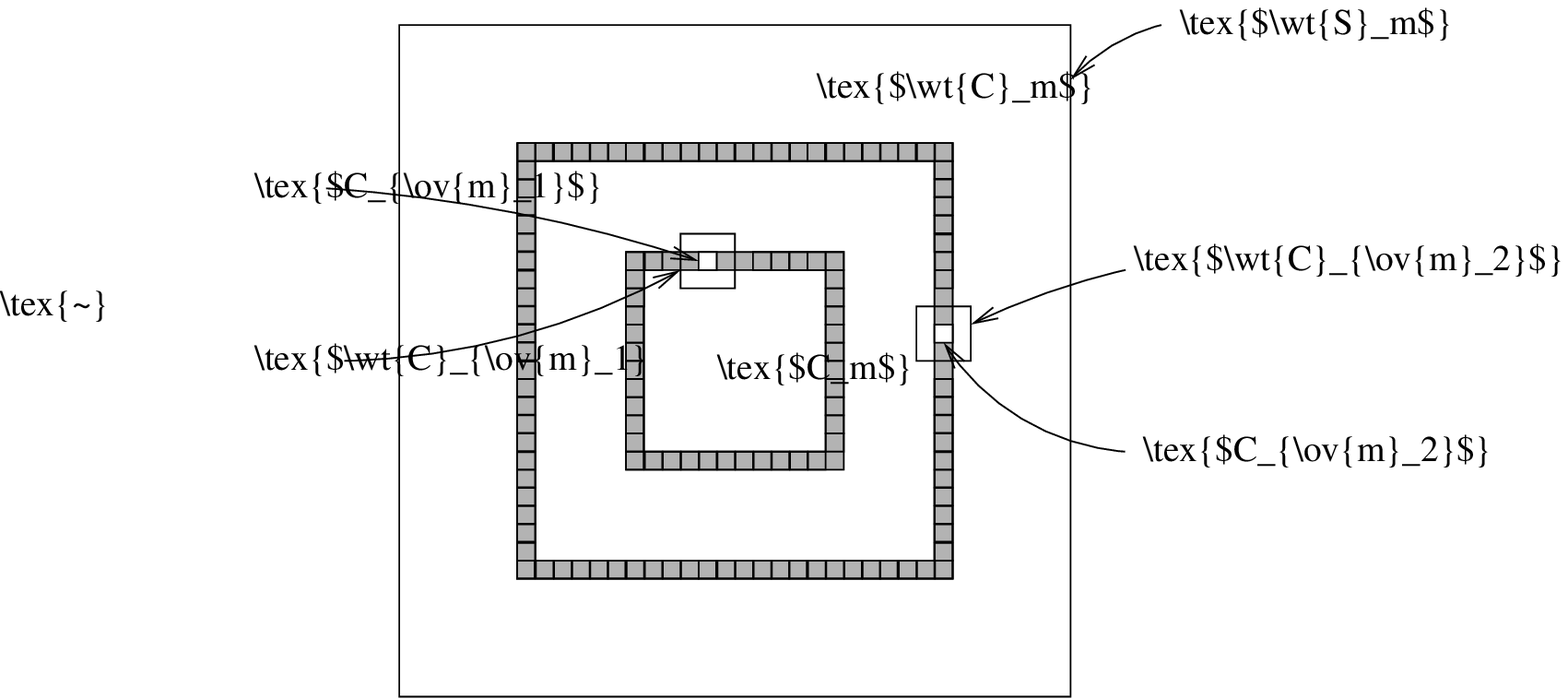}
\end{center}
\begin{center}
Fig.~1: An illustration of the boxes $C_{\ov{m}_i}$ and $\wt{C}_{\ov{m}_i}$, $i=1,2$.
\end{center}

\n
Given $m \in I_n$, $n \ge 0$, we consider $\Lambda_m$ the collection of subsets $\cT$ of $\bigcup_{0 \le k \le n} I_k$, such that setting $\cT^k = \cT \cap I_k$ one has
\begin{align}
& \cT^n = \{m\},\label{2.7}
\\[1ex]
&\mbox{any $m^\prime \in \cT^k$, $1 \le k \le n$ has two ``descendants'' $\ov{m}_1(m^\prime) \in \cH_1(m^\prime)$ and} \label{2.8}
\\
& \ov{m}_2(m^\prime) \in \cH_2(m^\prime), \;\mbox{such that $\cT^{k-1} = \bigcup_{m^\prime \in \cT^k} \{\ov{m}_1(m^\prime), \ov{m}_2(m^\prime)\}$}\,. \nonumber
\end{align}

\medskip
In other words any $\cT \in \Lambda_m$ has a natural structure of binary tree of depth $n$ with root $m (\in I_n)$ and for $0 \le k \le n$, $|\cT^k| = 2^{n-k}$, moreover the $\wt{C}_{m^\prime}$, for $m^\prime \in \cT^k$, are pairwise disjoint.

\medskip
Given $n \ge 0$, $m \in I_n$, and $\cT \in \Lambda_m$, one can attach to each $m^\prime \in \cT_m$ a binary tree $\cT_{m^\prime} \in \Lambda_{m^\prime}$, which, roughly speaking, consists of the descendants of $m^\prime$ in $\cT$:
\begin{equation}\label{2.9}
\cT_{m^\prime} = \{m^{\prime\prime} \in \cT; \;\wt{C}_{m^{\prime\prime}} \subseteq \wt{C}_{m^\prime}\},
\end{equation}
and for any $1 \le k \le n$, $m^\prime \in \cT^k$, one has the identity
\begin{equation}\label{2.10}
\cT_{m^\prime} = \{m^\prime\} \cup \cT_{\ov{m}_1(m^\prime)} \cup \cT_{\ov{m}_2(m^\prime)}, \; \mbox{(disjoint union)},
\end{equation}
where $\ov{m}_i(m^\prime) \in \cH_i(m^\prime)$, $i = 1,2$, are as in (\ref{2.8}), see also (\ref{2.6}).

\medskip
For $m \in I_n$ one has the following rough bound on the cardinality of $\Lambda_m$, the collection of binary trees attached to $m$:
\begin{equation}\label{2.11}
|\Lambda_m| \le (c\,\ell_0^{d-1})^2 (c \,\ell_0^{d-1})^4 \dots (c\,\ell_0^{d-1})^{2^n} = (c\,\ell_0^{d-1})^{2(2^n - 1)} = (c_0 \,\ell_0^{2(d-1)})^{2^n - 1} \,.
\end{equation}
We then introduce the events
\begin{equation}\label{2.12}
A^u_m = \{C_m \stackrel{\cV^u}{\longleftrightarrow} \wt{S}_m\}, \;\mbox{for $u \ge 0, m \in I_n, n \ge 0$}\,,
\end{equation}

\n
where the notation is similar as in (\ref{0.2}). The role of $\Lambda_m$ as a way to separate combinatorial complexity and probabilistic estimates comes in the next simple

\begin{lemma}\label{lem2.1} $(n \ge 0, m \in I_n, u \ge 0)$
\begin{equation}\label{2.13}
\begin{array}{l}
\IP [A^u_m] \le |\Lambda_m| \,p_n(u), \; \mbox{where}
\\[1ex]
p_n(u) = \sup\limits_{\cT \in \Lambda_m} \IP [A^u_{\cT}] \;\mbox{and} \;A^u_\cT = \bigcap\limits_{m^\prime \in \cT^0} A^u_{m^\prime}, \;\mbox{(recall $\cT^0 = \cT \cap I_0$)}\,.
\end{array}
\end{equation}
\end{lemma}

\begin{proof}
Observe that when $n \ge 1$, $m \in I_m$, any path in $\cV^u$ originating in $C_m$ and ending in $\wt{S}_m$ must go through some $C_{\ov{m}_1}, \ov{m}_1 \in \cH_1(m)$, reach $\wt{S}_{\ov{m}_1}$ and then go through some $C_{\ov{m}_2}, \ov{m}_2 \in \cH_2(m)$ and reach $\wt{S}_{\ov{m}_2}$. Hence one has the inclusion
\begin{equation*}
A^u_m \subseteq \textstyle{\bigcup\limits_{\ov{m}_i \in \cH_i(m), i=1,2}} A^u_{\ov{m}_1} \cap A^u_{\ov{m}_2} \,.
\end{equation*}

\n
Note that when $m \in I_0$, $A^u_m = A^u_\cT$, with $\cT = \{m\}$ the unique element of $\Lambda_m$.  Therefore with a straightforward induction on $n$, using the above inclusion, one finds that $A^u_m \subseteq \bigcup_{\cT \in \Lambda_m} A^u_\cT$, for all $m \in I_n$, $n \ge 0, u \ge 0$. The claim (\ref{2.13}) readily follows. 
\end{proof}

Unlike what was done in Section 3 of \cite{Szni07a} we will not estimate $\IP[A^u_m]$, $m \in I_n$, by induction on $n$, (this quantity does not depend on which $m \in I_n$ we consider, due to translation invariance). Instead we will estimate $p_n(u)$ by induction on $n$. This will lead to a finer specification of the possible long range dependence effects. As in \cite{Szni07a}, we will use sprinkling, i.e. we will control $p_{n+1}(u_{n+1})$ in terms of $p_n(u_n)$, along an increasing sequence $u_n$. The additional paths entering the random interlacement corresponding to the increase of $u_n$ into $u_{n+1}$ will enable to dominate long range dependence. The main step comes in the next proposition (compare with Proposition 3.1 of \cite{Szni07a}).

\begin{proposition}\label{prop2.2} $(d \ge 3)$

\medskip
There exist positive constants $c_1,c_2,c$ such that for $\ell_0 \ge c$, for increasing sequences $u_n, n \ge 0$, in $(0,\infty)$ and non-decreasing sequences $r_n, n \ge 0$, of positive integers, such that
\begin{equation}\label{2.14}
u_{n+1} \ge u_n (1 + c_1\,2^n \,\ell_0^{-(n+1)(d-2)})^{r_n + 1}, \;\mbox{for $n \ge 0$}\,,
\end{equation}
one has for all $n \ge 0$,
\begin{equation}\label{2.15}
p_{n+1} (u_{n+1}) \le p_n(u_{n+1}) \big(p_n(u_n) + u_n \,L_0^{d-2} (c_2^{n+1} \,\ell_0^{-(n+1)(d-2)}\big)^{r_n})\,.
\end{equation}

\medskip\n
(Note that $p_n (\cdot)$ is a non-increasing function and $p_n(u_{n+1}) \le p_n(u_n))$.
\end{proposition}

\begin{proof}
We consider some $n \ge 0$, $m \in I_{n+1}$, $\cT \in \Lambda_m$, and write $\ov{m}_1, \ov{m}_2$ for the unique elements of $\cH_1(m)$, $\cH_2(m)$ in $\cT^n (= \cT \cap I_n)$. We also consider $0 < u^\prime < u$, respectively playing the role of $u_n$ and $u_{n+1}$, as well as an integer $r \ge 1$, playing the role of $r_n$. We write $\ov{\cT}_1$ and $\ov{\cT}_2$ in place of $\cT_{\ov{m}_1}$ and $\cT_{\ov{m}_2}$. We also define
\begin{equation}\label{2.16}
V = \wh{C}_1 \cup \wh{C}_2, \;\mbox{where} \; \wh{C}_i = \textstyle{\bigcup\limits_{m^\prime \in \ov{\cT}_i \cap I_0}} \wt{C}_{m^\prime} \subseteq \wt{C}_{\ov{m}_i}, \;\mbox{for $i = 1,2$}\,.
\end{equation}
We then introduce the decomposition, see (\ref{1.7}) for the notation,
\begin{equation}\label{2.17}
\mu_{V,u} = \mu_{1,1} + \mu_{1,2} + \mu_{2,1} + \mu_{2,2},
\end{equation}
where for $i \not= j$ in $\{1,2\}$ we have set
\begin{align*}
\mu_{i,j} & = 1\{X_0 \in \wh{C}_i, \; H_{\wh{C}_j} < \infty\} \,\mu_{V,u}, \;\mbox{and}
\\
\mu_{i,i} & = 1\{X_0 \in \wh{C}_i, \; H_{\wh{C}_j} = \infty\} \,\mu_{V,u}.
\end{align*}

\n
This is similar to (3.14) of \cite{Szni07a}, except that $\wh{C}_i$, $i = 1,2$, now plays the role of $\wt{C}_{\ov{m}_i}$, $i = 1,2$ in (3.14) of \cite{Szni07a}. Similarly with $\mu_{V,u^\prime}$, and $\mu_{V, u^\prime,u}$, see (\ref{1.7}) for the notation, we define with analogous formulas
\begin{align}
&\;\;  \mu_{V,u^\prime} = \mu^\prime_{1,2} + \mu^\prime_{1,2} + \mu^\prime_{2,1} + \mu^\prime_{2,2} \,,\nonumber
\\[0.5ex]
&\mu_{V,u^\prime,u} = \mu^*_{1,1} + \mu^*_{1,2} + \mu^*_{2,1} + \mu^*_{2,2}, \;\mbox{so that}\nonumber
\\[0.5ex]
&\mbox{$\mu^\prime_{i,j}$, $\mu^*_{i,j}$, $1 \le i, j \le 2$, are independent Poisson point processes on $W_+$}\,.\label{2.18}
\end{align}

\medskip\n
When $\mu$ is a random point process on $W_+$ defined on $\Omega$, (i.e.~a measurable map from $\Omega$ into the space of point measures on $W_+)$, it will be convenient to define for $\ov{\cT} \in \Lambda_{\ov{m}}$, with $\ov{m} \in I_n$
\begin{equation}\label{2.19}
\begin{split}
A_{\ov{\cT}}(\mu) = \textstyle{\bigcap\limits_{m^\prime \in \ov{\cT} \cap I_0}} \big\{ \o \in \Omega ; &\; \mbox{there is a path in $\wt{C}_{m^\prime} \backslash \big(\bigcup\limits_{w \in \, {\rm Supp}\,\mu(\omega)} w(\IN)\big)$} 
\\[-3ex]
&\; \mbox{joining $C_{m^\prime}$ with $\wt{S}_{m^\prime}\big\}$},
\end{split}
\end{equation}
In particular for $\ov{m} \in I_n$, $\ov{\cT} \in \Lambda_{\ov{m}}$, one finds that, cf.~(\ref{1.12}),
\begin{equation*}
A^u_{\ov{\cT}} = A_{\ov{\cT}} (\mu_{K,u}), \;\mbox{for any finite} \;K \supseteq \textstyle{\bigcup\limits_{m^\prime \in \ov{\cT} \cap I_0} } \wt{C}_{m^\prime} \,.
\end{equation*}

\n
With (\ref{2.10}) it follows that $A^u_{\cT} = A^u_{\ov{\cT}_1} \cap A^u_{\ov{\cT}_2}$, and taking into account that $w \in {\rm Supp} \,\mu_{2,2}$ implies that $w(\IN) \cap \wh{C}_1 = \emptyset$, we have
\begin{align*}
A^u_{\cT} & = A^u_{\ov{\cT}_1} \cap A^u_{\ov{\cT}_2} = A^u_{\ov{\cT}_1}  (\mu_{1,1}  + \mu_{1,2} + \mu_{2,1}) \cap A^u_{\ov{\cT}_2} (\mu_{V,u})
\\[1ex]
& \subseteq A^u_{\ov{\cT}_1} (\mu_{1,1} + \mu_{1,2} + \mu_{2,1}) \cap A^u_{\ov{\cT}_2} (\mu_{2,2})\,.
\end{align*}

\n
Using the independence of the $\mu_{i,j}$, $1 \le i,j \le 2$, we find that
\begin{equation}\label{2.20}
\begin{split}
\IP [A^u_\cT] & \le \IP [A_{\ov{\cT}_1} (\mu_{1,1} + \mu_{1,2} + \mu_{2,1})] \;\IP [A_{\ov{\cT}_2} (\mu_{2,2})\big]
\\[1ex]
& = \IP [A^u_{\ov{\cT}_1} ] \;\IP [A_{\ov{\cT}_2} (\mu_{2,2})] \le p_n(u) \;\IP[A_{\ov{\cT}_2}  (\mu_{2,2})]\,.
\end{split}
\end{equation}

\n
We will now bound $\IP [A_{\ov{\cT}_2} (\mu_{2,2})] \stackrel{(\ref{1.8})}{=} \IP[A_{\ov{\cT}_2} (\mu^\prime_{2,2} + \mu^*_{2,2})]$ in terms of $p_n(u^\prime)$ when $u - u^\prime$ is substantial enough. For this purpose  cf.~(\ref{2.35}) below, we will dominate the influence on $\wh{C}_2$ of $\mu^\prime_{2,1} + \mu^\prime_{2,1}$ by $\mu^*_{2,2}$. This has the same spirit as what appears in Proposition 3.1 of \cite{Szni07a}, however $\wh{C}_2$ now replaces $\wt{C}_{\ov{m}_2}$. The sprinkling technique will come into action during this step.  

\medskip
We introduce the $\ell^\infty$-neighborhood of size $\frac{L_{n+1}}{10}$ of $\wt{C}_{\ov{m}_2}$ ($\frac{L_{n+1}}{10}$ is an integer since $\ell_0 \ge 100$ is a multiple of $10$, cf.~(\ref{2.1})):
\begin{equation}\label{2.21}
U = \Big\{z \in \IZ^d; d (z, \wt{C}_{\ov{m}_2}) \le \mbox{\f $\dis\frac{L_{n+1}}{10}$}\Big\} \,.
\end{equation}
We then consider the times of successive returns to $\wh{C}_2$ and departures from $U$:
\begin{equation}\label{2.22}
\begin{split}
R_1 & = H_{\wh{C}_2}, \; D_1 = T_U \circ \theta_{R_1} + R_1, \;\mbox{and by induction}
\\[0.5ex]
R_{k+1} & = R_1 \circ \theta_{D_k} + D_k, \;D_{k+1} = D_1 \circ \theta_{D_k} + D_k, \;\mbox{for $k \ge 1$},
\end{split}
\end{equation}
so that $0 \le R_1 \le D_1 \le \dots \le R_k \le D_k \le \dots \le \infty$.

\medskip
We then introduce the decompositions:
\begin{equation}\label{2.23}
\begin{split}
\mu^\prime_{2,1} & = \dsl_{1 \le \ell \le r} \;\rho^\ell_{2,1} + \ov{\rho}_{2,1}, \;\mu^\prime_{1,2} = \dsl_{1 \le \ell \le r} \,\rho^\ell_{1,2} + \ov{\rho}_{1,2} \,,
\\[0.5ex]
\mu^*_{2,2} & = \dsl_{1 \le \ell \le r} \,\rho^\ell_{2,2} + \ov{\rho}_{2,2}\,,
\end{split}
\end{equation}
where for $i \not= j$ in $\{1,2\}$ and $\ell \ge 1$, we have set
\begin{align*}
\rho^\ell_{i,j} & = 1\{R_\ell < D_\ell < R_{\ell + 1}= \infty\}\;\mu^\prime_{i,j}, \quad \! \ov{\rho}_{i,j} = 1\{R_{r+1} < \infty\}\,\mu^\prime_{i,j}, \;\mbox{and}
\\
\rho^\ell_{2,2} & = 1\{R_\ell < D_\ell < R_{\ell + 1}= \infty\}\;\mu^*_{2,2}, \;\;\ov{\rho}_{2,2} = 1\{R_{r+1} < \infty\}\,\mu^*_{2,2}.
\end{align*}

\n
As a result of (\ref{2.18}) and the above formulas we see that:
\begin{equation}\label{2.24}
\begin{array}{l}
\mbox{$\mu^\prime_{2,2}, \rho^\ell_{i,j}, 1 \le \ell \le r, \ov{\rho}_{i,j}, 1 \le i,j \le 2$, with $i$ or $j \not= 1$, are independent}
\\
\mbox{Poisson point processes on $W_+$}\,.
\end{array}
\end{equation}

\medskip\n
We denote with $\ov{\xi}_{2,1}$ and $\ov{\xi}_{1,2}$ the respective intensity measures of $\ov{\rho}_{2,1}$ and $\ov{\rho}_{1,2}$. We have
\begin{equation}\label{2.25}
\begin{split}
\ov{\rho}_{2,1}(W_+) &= u^\prime \,P_{e_V} [X_0 \in \wh{C}_2, \;H_{\wh{C}_1} < \infty, \;R_{r+1} < \infty]
\\
& \le   u^\prime \;{\rm cap} (\wh{C}_2) \sup\limits_{x \in \wh{C}_2} \;P_x [R_{r+1} < \infty]   \overset{\stackrel{\mbox{\f strong}}{\mbox{\f Markov}}}{\le}  \!\!\! u^\prime \;{\rm cap}(\wh{C}_2) (\sup\limits_{x \in U^c} \;P_x [H_{\wh{C}_2} < \infty])^r.
\end{split}
\end{equation}

\n
Note that with standard estimates on the Green function, cf.~\cite{Lawl91}, p.~31, as well as (\ref{1.3}), (\ref{1.4}), (\ref{1.6}), we find that
\begin{equation}\label{2.26}
\sup\limits_{x \in U^c} \;P_x [H_{\wh{C}_2} < \infty] \le c\,2^n \;\dis\frac{L^{d-2}_0}{L_{n+1}^{d-2}} \stackrel{(\ref{2.1})}{=} c\,2^n \,\ell_0^{-(n+1)(d-2)}\,.
\end{equation}
Therefore with (\ref{2.25}), and (\ref{1.3}), (\ref{1.6}), we find that
\begin{equation}\label{2.27}
\begin{split}
\ov{\xi}_{2,1}(W_+) & \le c\,u^\prime \;2^n \,L_0^{d-2} \Big( c\,2^n \,\ell_0^{-(n+1)(d-2)}\Big)^r
\\[1ex]
&\!\! \overset{r \ge 1}{\le}  u^\prime \,L_0^{d-2} \Big(c\,4^n \,\ell_0^{-(n+1)(d-2)}\Big)^r\,.
\end{split}
\end{equation}
In a similar fashion, we also have
\begin{equation}\label{2.28}
\begin{split}
\ov{\xi}_{1,2}(W_+) & = u^\prime \;P_{e_V} [X_0 \in \wh{C}_1, H_{\wh{C}_2} < \infty, \;R_{r+1} < \infty]
\\
& \le u^\prime \,L_0^{d-2} \Big(c\,4^n \,\ell_0^{-(n+1)(d-2)}\Big)^r \,.
\end{split}
\end{equation}

\n
We will now prove that the trace on $\wh{C}_2$ of paths in the support of $\sum_{1 \le \ell \le r} \rho_{2,1}^\ell$ and $\sum_{1 \le \ell \le r} \rho_{1,2}^\ell$ is dominated by the corresponding trace of paths in the support of $\mu^*_{2,2}$ if $u-u^\prime$ is not too small.

\medskip
We consider $W_f$ the collection of finite paths in $\IZ^d$ and for $\ell \ge 1$, the measurable map $\phi^\ell$ from $\{D_\ell < R_{\ell + 1} = \infty\} \subseteq W_+$ into $W_f^{\times \ell}$ such that
\begin{equation}\label{2.29}
\phi^\ell(w) = (w(R_k + \cdot)_{0 \le \cdot \le D_k - R_k})_{1 \le k \le \ell} \in W_f^{\times \ell}, \;\mbox{for} \; w \in \{D_\ell < R_{\ell + 1} = \infty\}\,.
\end{equation}

\medskip\n
In other words $\phi^\ell(w)$ keeps track of the parts of the trajectory $w$ going from the successive returns to $\wh{C}_2$ up to departure from $U$. We view the $\rho^\ell_{i,j}$, $i$ or $j \not= 1$, with $\ell \ge 1$ fixed, as point processes on $\{D_\ell < R_{\ell + 1} = \infty\} \subseteq W_+$, and denote with $\wt{\rho}\,{^\ell_{i,j}}$ the respective images under $\phi^\ell$, which are Poisson point processes on $W_f^{\times \ell}$. We write $\wt{\xi}\,{_{i,j}^\ell}$ for the corresponding intensity measures. With (\ref{2.24}) it follows that
\begin{equation}\label{2.30}
\begin{array}{l}
\mu^\prime_{2,2}, \; \wt{\rho}\,{^\ell_{i,j}}, \;1 \le \ell \le r, \;\ov{\rho}_{i,j}, \;1 \le i, j \le 2, 
\\[1ex]
\mbox{$i$ or $j \not= 1$, are independent point processes}.
\end{array}
\end{equation}

\n
Our next step is the following lemma, which is an adaption of Lemma 3.2 of \cite{Szni07a}.

\bigskip\n
\begin{lemma}\label{lem2.3}
For $\ell_0 \ge c$, for all $n \ge 0$, $m \in I_{n+1}$, $\cT \in \Lambda_m$, $x \in \partial U$, $y \in \partial_{\rm int} \,\wh{C}_2$, one has
\begin{align}
&P_x [H_{\wh{C}_1} < R_1 < \infty, X_{R_1} = y] \le c\,2^n \,\ell_0^{-(n+1)(d-2)} P_x [H_{\wh{C}_1} > R_1, X_{R_1} = y],\label{2.31}
\\[1ex]
&P_x [H_{\wh{C}_1}  < \infty, R_1 = \infty] \le c\,2^n \,\ell_0^{-(n+1)(d-2)} P_x [R_1= \infty = H_{\wh{C}_1}].\label{2.32}
\end{align}
\end{lemma}

\begin{proof}
We begin with the proof of (\ref{2.31}). For $z \in \partial U$, $y \in \partial_{\rm int} \wh{C}_2$, one finds with the strong Markov property that
\begin{align*}
&P_z [H_{\wh{C}_1} < R_1 < \infty, X_{R_1} = y]  = E_z\big[H_{\wh{C}_1} < R_1 , P_{X_{H_{\wh{C}_1}}} [R_1 < \infty, X_{R_1} = y]\big] =
\\
&P_z \big[H_{\wh{C}_1}  <  R_1, E_{X_{H_{\wh{C}_1}}} \big[H_{\partial U} < \infty, P_{X_{H_{\partial U}}} [R_1 < \infty, X_{R_1} = y]\big]\big],
\end{align*}

\n
where we have used in the last step that for $z^\prime \in \wh{C}_1$, $P_{z^\prime}$-a.s., $R_1 = H_{\partial U} + R_1 \circ \theta_{H_{\partial U}}$. As a result we obtain
\begin{equation}\label{2.33}
\begin{array}{l}
\sup\limits_{z \in \partial U} P_z[H_{\wh{C}_1} < R_1 < \infty, X_{R_1} = y] \le
\\[2ex]
 \sup\limits_{z \in \partial U} P_z [H_{\wh{C}_1} < \infty] \;\sup\limits_{z \in \partial U} P_z [R_1 < \infty, X_{R_1} = y] \le
\\[2ex]
c \,2^n \,\ell_0^{-(n+1)(d-2)} \sup\limits_{z \in \partial U} \;P_z [R_1 < \infty, X_{R_1} = y] ,
\end{array}
\end{equation}
with a similar bound as in (\ref{2.26}) in the last step.

\medskip
Now observe that $P_z [R_1 < \infty, X_{R_1}=y] = P_z [H_{\wh{C}_2} < \infty, X_{H_{\wh{C}_2}} = y]$, $z \in \wh{C}_2^c$ is a positive harmonic function and with Harnack inequality, cf.~Theorem 1.7.2, p.~42 of \cite{Lawl91}, together with a standard covering argument we have
\begin{equation*}
\sup\limits_{z \in \partial U} P_z[R_1 < \infty, X_{R_1} = y] \le c \inf\limits_{z \in \partial U} P_z[R_1 < \infty, X_{R_1} = y] \,.
\end{equation*}
Thus coming back to (\ref{2.33}) we find that
\begin{equation*}
\begin{array}{l}
\sup\limits_{z \in \partial U} P_z[H_{\wh{C}_1} < R_1 < \infty, X_{R_1} = y] \le c^\prime \,2^n\,\ell_0^{-(n+1)(d-2)} \inf\limits_{z \in \partial U} P_z[R_1 < \infty, X_{R_1} = y] =
\\[2ex]
c^\prime \,2^n\,\ell_0^{-(n+1)(d-2)} \inf\limits_{z \in \partial U} \;\big(P_z [H_{\wh{C}_1} < R_1 < \infty, X_{R_1} = y] + P_z [H_{\wh{C}_1} > R_1, X_{R_1} = y]\big)\,.
\end{array}
\end{equation*}

\n
For large $\ell_0$, we find that $c^\prime \,2^n\,\ell_0^{-(n+1)(d-2)} \le \frac{1}{2}$, for all $n \ge 0$, with $c^\prime$ as above and hence for $x \in \partial U$:
\begin{equation*}
P_x[H_{\wh{C}_1} < R_1 < \infty, X_{R_1} = y] \le 2c^\prime \, 2^n\,\ell_0^{-(n+1)(d-2)} P_x [H_{\wh{C}_1} > R_1, X_{R_1} = y],
\end{equation*}
and this proves (\ref{2.31}).

\medskip
We then turn to the proof of (\ref{2.32}) which is more straightforward. We observe that
\begin{equation*}
\inf\limits_{x \in \partial U} P_x[R_1 = \infty, H_{\wh{C}_1} = \infty] = \inf\limits_{x \in \partial U} P_x [H_V = \infty] \ge c\,,
\end{equation*}

\n
using the invariance principle to let the walk move at a distance of $\wt{C}_{\ov{m}_1} \cup \wt{C}_{\ov{m}_2}$, which is a multiple of $L_{n+1}$, as well as (\ref{1.5}), (\ref{1.6}) and standard bounds on the Green function. On the other hand the left-hand side of (\ref{2.32}) with a similar bound as in (\ref{2.26}) is smaller than $c\,2^n \,\ell_0^{-(n+1)(d-2)}$ and the claim follows.
\end{proof}

The main control on the intensity measure $\wt{\xi}^\ell_{1,2} + \wt{\xi}^\ell_{2,1}$ of $\wt{\rho}\,^\ell_{1,2} + \wt{\rho}\,^\ell_{2,1}$ in terms of the intensity measure $\wt{\xi}^\ell_{2,2}$ of $\wt{\rho}\,^\ell_{2,2}$ comes from

\begin{lemma}\label{lem2.4} $(\ell_0 \ge c)$
\begin{equation}\label{2.34}
\wt{\xi}^\ell_{1,2} + \wt{\xi}^\ell_{2,1} \le \mbox{\f $\dis\frac{u^\prime}{u-u^\prime}$}
 \Big[\big(1 + c_1 \,2^n\,\ell_0^{-(n+1)(d-2)}\big)^{\ell +1} - 1\Big] \wt{\xi}^\ell_{2,2}, \;\mbox{for $\ell \ge 1$} \,.
\end{equation}
\end{lemma}

\begin{proof}
This is a repetition of the proof of Lemma 3.3 of \cite{Szni07a}, with $\wh{C}_i$, $i = 1,2$, replacing $\wt{C}_{\ov{m}_i}$, $i = 1,2$ in \cite{Szni07a}, and $c\,2^n \,\ell_0^{-(n+1)(d-2)}$ replacing $c\,\ell_n^{-(d-2)}$ in the notation of \cite{Szni07a}, thanks to (\ref{2.31}), (\ref{2.32}) of Lemma \ref{lem2.3} above.
\end{proof}

We now suppose that $\ell_0 \ge c$, so that the Lemmas \ref{lem2.3} and \ref{lem2.4} apply, and that
\begin{equation}\label{2.35}
\begin{array}{l}
u \ge \Big(1 + c_1 \,2^n\,\ell_0^{-(n+1)(d-2)}\Big)^{r+1} u^\prime, 
\\
\Big(\mbox{hence} \;\mbox{\f $\dis\frac{u^\prime}{u-u^\prime}$} \Big[\Big(1 + c_1 \, 2^n\,\ell_0^{-(n+1)(d-2)}\Big)^{r+1} - 1\Big] \le 1\Big)\,.
\end{array}
\end{equation}

\n
In our present notation, (\ref{2.35}) coincides with (\ref{2.14}). We now bound $\IP[A_{\ov{\cT}_2}(\mu_{2,2})]$ in terms of $p_n (u^\prime) \ge \IP[A^{u^\prime}_{\ov{\cT}_2}]$ as follows. We first express the trace of $\cI^{u^\prime}$ on $\wh{C}_2$ as the union
\begin{equation} \label{2.36}
\cI^{u^\prime}  \cap \wh{C}_2 = \cI^\prime \cup \wt{\cI} \cup \ov{\cI}, \;\mbox{where}
\end{equation}

\medskip
\begin{equation} \label{2.37}
\begin{split}
\cI^\prime & = \textstyle{\bigcup\limits_{w \in \,{\rm Supp} (\mu^\prime_{2,2})}} w(\IN) \cap \wh{C}_2\,, 
\\
\wt{\cI} & = \textstyle{\bigcup\limits_{1 \le \ell \le r}}\quad   \textstyle{\bigcup\limits_{(w_1,\dots,w_\ell) \in \,{\rm Supp}\, \wt{\rho}^\ell_{1,2} + \wt{\rho}^\ell_{2,1}}} ({\rm range} \;w_1 \cup \dots \cup \;{\rm range} \;w_\ell) \cap \wh{C}_2\,,
\\
\ov{\cI} & = \textstyle{\bigcup\limits_{w \in \,{\rm Supp}\, \ov{\rho}_{1,2} + \ov{\rho}_{2,1}}}  w(\IN) \cap \wh{C}_2 \,. 
\end{split}
\end{equation}

\n
If we now define $\cI^*$ replacing $\wt{\rho}\,^\ell_{1,2} + \wt{\rho}\,^\ell_{2,1}$ by $\wt{\rho}\,^\ell_{2,2}$ in the second line of (\ref{2.37}), we find that with (\ref{2.30})
\begin{equation}\label{2.38}
\mbox{the random sets $\cI^\prime, \wt{\cI}, \ov{\cI}, \cI^*$ are independent under $\IP$} \,.
\end{equation}
Now with (\ref{2.34}), (\ref{2.35}) we also find that
\begin{equation}\label{2.39}
\mbox{$\wt{\cI}$ is stochastically dominated by $\cI^*$}\,.
\end{equation}
As as result we can write that
\begin{equation}\label{2.40}
\begin{array}{l}
\IP \Big[A_{\ov{\cT}_2} \Big(\mu^\prime_{2,2} + \dsl_{1 \le \ell \le r} \;\rho^\ell_{2,2}\Big)\Big] =
\\[3ex]
 \IP\big[\mbox{for each $m^\prime \in \ov{\cT}_2 \cap I_0$, there is a path in $\wt{C}_{m^\prime} \backslash (\cI^\prime \cup \cI^*)$ from $C_{m^\prime}$ to $\wt{S}_{m^\prime}]$} \le
\\[2ex]
\IP\big[\mbox{for each $m^\prime \in \ov{\cT}_2 \cap I_0$, there is a path in $\wt{C}_{m^\prime} \backslash (\cI^\prime \cup\wt{\cI})$ from $C_{m^\prime}$ to $\wt{S}_{m^\prime} \big]$} =
\\[1ex]
\IP\Big[A_{\ov{\cT}_2} \Big(\mu^\prime_{2,2} + \dsl_{1 \le \ell \le r} \;\rho^\ell_{2,1} + 
\rho^\ell_{1,2}\Big)\Big] ,
\end{array}
\end{equation}
and hence that
\begin{equation}\label{2.41}
\begin{array}{l}
\hspace{-3ex} \IP [A_{\ov{\cT}_2} (\mu_{2,2})] = \IP [A_{\ov{\cT}_2}(\mu^\prime_{2,2} + \mu^*_{2,2})] \le \IP \Big[A_{\ov{\cT}_2} \Big(\mu^\prime_{2,2} + \dsl_{1 \le \ell \le r} \;\rho^\ell_{2,2}\Big)\Big] \stackrel{(\ref{2.40})}{\le} 
\\[2ex]
\hspace{-3ex}   \IP \Big[A_{\ov{\cT}_2} \Big(\mu^\prime_{2,2} + \dsl_{1 \le \ell \le r} \;\rho^\ell_{2,1} + \rho^\ell_{1,2}\Big)\Big] \le
\\[1ex]
\hspace{-3ex} \IP [A_{\ov{\cT}_2} (\mu^\prime_{2,2} + \mu^\prime_{2,1} + \mu^\prime_{1,2}), \ov{\rho}_{2,1} = \ov{\rho}_{1,2} = 0] + \IP[\ov{\rho}_{2,1} \;\mbox{or} \;\ov{\rho}_{1,2} \not= 0] = 
\\[1ex]
\hspace{-3ex}  \IP [A_{\ov{\cT}_2} (\mu_{V,u^\prime}), \ov{\rho}_{2,1} = \ov{\rho}_{1,2} = 0] + \IP [\ov{\rho}_{2,1} \;\mbox{or} \; \ov{\rho}_{1,2} \not= 0] \stackrel{{\rm below} \;(\ref{2.19})}{\le} 
\\[1ex]
 \hspace{-3ex}  p_n(u^\prime) + \ov{\xi}_{2,1}(W_+) + \ov{\xi}_{1,2}(W_+) \stackrel{(\ref{2.27}),(\ref{2.28})}{\le} p_n(u^\prime) + 2u^\prime \,L_0^{d-2} (c\,4^n\,\ell_0^{-(n+1)(d-2)})^r\,.
\end{array}
\end{equation}
This proves (\ref{2.15}). 
\end{proof}

\n
We will now consider sequences $u_n, n\ge 0$, and $r_n, n\ge 0$, of the form
\begin{align}
u_n =&\; u_0 \,\exp\Big\{c_1 \dsl_{0 \le k < n} (r_k + 1) \,2^k\,\ell_0^{-(k+1)(d-2)}\Big\}\,, \label{2.42}
\\
r_n = &\;  r_0\,2^n, \label{2.43}
\end{align}

\n
where $u_0$ is a positive number and $r_0$ a positive integer. The choice (\ref{2.42}) ensures in particular that (\ref{2.14}) holds. Observe also that the increasing sequence $u_n$ has a finite limit
\begin{equation}\label{2.44}
u_\infty = u_0 \exp\Big\{\mbox{\f $\dis\frac{c_1}{\ell_0^{d-2}}$}\;\Big(\mbox{\f $\dis\frac{r_0}{1-4 \,\ell_0^{-(d-2)}}$} + \mbox{\f $\dis\frac{1}{1-2 \ell_0^{-(d-2)}}$}\Big)\Big\}\,,
\end{equation}

\n
where we recall that $\ell_0$ is at least $100$ and even.

\medskip
The next proposition encapsulates bounds on $p_n(u_n)$, which can be propagated with the help of Proposition \ref{prop2.2}, if we can initiate the induction. This will be our main tool in the proof of Theorem \ref{theo0.1} in the next section.

\medskip
With (\ref{2.44}) we will view $u_\infty$ as a function of $u_0,r_0, \ell_0$.

\begin{proposition}\label{prop2.5}
There exists a positive constant $c$ such that when $u_0 > 0$, $r_0 \ge 1$, $\ell_0 \ge c$, $L_0 \ge 1$, $K_0 > \log 2$ are such that in the notation of {\rm (\ref{2.15}), (\ref{2.44})},
\begin{align}
&u_\infty \,L_0^{d-2} \vee e^{K_0} \le \Big(\mbox{\f $\dis\frac{\ell_0^{d-2}}{c_2}$}\Big)^{\frac{r_0}{2}}, \;\mbox{and} \label{2.45}
\\[1ex]
&p_0 (u_0) \le e^{-K_0}, \label{2.46}
\intertext{then}
&p_n(u_n) \le e^{-(K_0 - \log 2)2^n}, \;\mbox{for each $n \ge 0$} \,. \label{2.47}
\end{align}
\end{proposition}

\begin{proof}
We assume $\ell_0 \ge c$, so that Proposition \ref{prop2.2} applies. Note that due to (\ref{2.45}) one has $\ell_0^{d-2} \ge c_2$.

\medskip
The inequality (\ref{2.47}) will result from an induction argument relying on (\ref{2.15}). For this purpose we observe that the last term in the right-hand side of (\ref{2.15}) can be bounded as follows:
\begin{equation*}
u_n\,L_0^{d-2} \big(c_2^{n+1} \,\ell_0^{-(n+1)(d-2)}\big)^{r_n} \le u_\infty \,L_0^{d-2} \Big(\mbox{\f $\dis\frac{c_2}{\ell_0^{d-2}}$}\Big)^{r_0} \big(c_2 \, \ell_0^{-(d-2)}\big)^{nr_n} \stackrel{(\ref{2.45})}{\le} \big(c_2 \, \ell_0^{-(d-2)}\big)^{\frac{r_n}{2}} \,.
\end{equation*}

\n
As a result we see that (\ref{2.15}) implies that for all $n \ge 0$,
\begin{equation}\label{2.48}
p_{n+1} (u_{n+1}) \le p_n(u_n) \big(p_n(u_n) + \big(c_2 \,\ell_0^{-(d-2)}\big)^{\frac{r_0}{2} 2^n}\big)\,.
\end{equation}

\n
We then define by induction a sequence $K_n$, $n \ge 0$, such that
\begin{equation}\label{2.49}
K_n = K_0 - \dsl_{0 \le n^\prime < n} \;\mbox{\f $\dis\frac{1}{2^{n^\prime + 1}}$} \log \big(1 + e^{K_{n^\prime} \,2^{n^\prime}}\big(c_2 \,\ell_0^{-(d-2)}\big)^{\frac{r_0}{2} 2^{n^\prime}}\big), \; \mbox{for $n \ge 1$}\,.
\end{equation}
Note that $K_n \le K_0$ so that
\begin{equation}\label{2.50}
\begin{split}
K_n & \ge K_0 - \dsl_{n^\prime \ge 0} 2^{-(n^\prime + 1)} \log \big(1 + e^{K_0 2^{n^\prime}} \big(c_2 \,\ell_0^{-(d-2)}\big)^{\frac{r_0}{2} 2^{n^\prime}}\big)
\\[1ex]
&\!\!\!\! \stackrel{(\ref{2.45})}{\ge} K_0 - \dsl_{n^\prime \ge 0} 2^{-(n^\prime + 1)} \log 2 = K_0 - \log 2 > 0\,.
\end{split}
\end{equation}
We will now check by induction that
\begin{equation}\label{2.51}
p_n(u_n) \le e^{-K_n 2^n}, \;\mbox{for all $n \ge 0$}\,.
\end{equation}

\n
In view of (\ref{2.50}) this will imply (\ref{2.47}). The assumption (\ref{2.46}) ensures that (\ref{2.51}) holds when $n = 0$. We then assume that it holds for $n$ and find with (\ref{2.48}) that
\begin{align*}
p_{n+1}(u_{n+1}) & \le e^{-K_n 2^n} (e^{-K_n 2^n} + \big(c_2 \,\ell_0^{-(d-2)}\big)^{\frac{r_0}{2} 2^n}\big)
\\[1ex]
& = e^{-K_n 2^{n+1}} \big(1 + e^{K_n 2^n} \big(c_2 \,\ell_0^{-(d-2)}\big)^{\frac{r_0}{2} 2^n}\big)
\\[1ex]
& = e^{-2^{n+1} [K_n - \frac{1}{2^{n+1}} \log (1 + e^{K_n 2^n} (c_2 \,\ell_0^{-(d-2)})^{\frac{r_0}{2} 2^n})]}\stackrel{(\ref{2.49})}{=} e^{-K_{n+1} 2^{n+1}} \,.
\end{align*}

\n
This proves (\ref{2.51}) for $n+1$, and concludes the proof by induction of (\ref{2.51}) and hence of Proposition \ref{prop2.5}.
\end{proof}

\section{Denouement}
\setcounter{equation}{0}

We will now prove Theorem \ref{theo0.1} with the help of Proposition \ref{prop2.5} of the previous section. We recall the notation from (\ref{0.2}).

\bigskip\n
{\it Proof of Theorem {\rm \ref{theo0.1}}:} We consider $u > u_{**}$ and define $u_0 = \frac{1}{2} (u + u_{**}) \in (u_{**},u)$. We know from (\ref{0.2}) that for some $\varepsilon(u) \in (0,1)$, one has
\begin{equation}\label{3.1} \;\
\lim\limits_{L \r \infty} L^\ve \IP[B(0,L) \stackrel{\cV^{u_0}}{\longleftrightarrow} S(0,2L)] = 0\,.
\end{equation}
This readily implies that in the notation of (\ref{2.4}), (\ref{2.12})
\begin{equation}\label{3.2}
\lim\limits_{L_0 \r \infty} L_0^\ve \IP[A^{u_0}_m] = 0, \;\mbox{for arbitrary $m \in I_0$} \,.
\end{equation}

\n
We select the parameters $r_0, K_0, \ell_0$, which appear in Proposition \ref{prop2.5} as follows:
\begin{align}
r_0 & = \Big[\mbox{\f $\dis\frac{12}{\ve}$} \;(d-1)\Big] + 1 \label{3.3}
\\[1ex]
K_0 & = \log \big(c_0 \,\ell_0^{2(d-1)}\big) + 2 \log 2, \;\;  \mbox{(see (\ref{2.11}) for the notation)} \label{3.4}
\\[1ex]
\ell_0 & = 200 \big(\big[L_0^{\frac{\ve}{3(d-1)}}\big] + 1\big) \,.\label{3.5}
\end{align}

\n
There remains to select $L_0$. We will now see that
\begin{equation}\label{3.6}
\mbox{when $L_0 \ge c(u)$,  (\ref{2.45}) and (\ref{2.46}) are fulfilled}\,.
\end{equation}

\n
We first observe that in view of (\ref{2.44}), for $L_0 \ge c(u)$, $u_\infty < u$. As a result for $L_0 \ge c^\prime(u)$ we find that
\begin{equation*}
\Big(\mbox{\f $\dis\frac{\ell_0^{d-2}}{c_2}$}\Big)^{\frac{r_0}{2}} \stackrel{(\ref{3.3}),(\ref{3.5})}{\ge} \Big(c \,L_0^{\frac{(d-2)\ve}{3(d-1)}}\Big)^{\frac{6}{\ve}\,(d-1)} = c(u) \,L_0^{2(d-2)} \ge u\,L_0^{d-2} \ge u_\infty\,L_0^{d-2},
\end{equation*}
and that
\begin{equation*}
e^{K_0} = 4 \,c_0 \,\ell_0^{2(d-1)} \le \Big(\mbox{\f $\dis\frac{\ell_0^{d-2}}{c_2}$}\Big)^{\frac{r_0}{2}} ,
\end{equation*}
when we have used that $r_0(d-2) > 4(d-1)$.

\medskip
This takes care of condition (\ref{2.45}). As for (\ref{2.46}) we note that for $L_0 \ge c(u)$, one has
\begin{equation*}
e^{-K_0} = (4 \,c_0 \,\ell_0^{2(d-1)})^{-1} \stackrel{(\ref{3.5})}{\ge} c\,L_0^{-\frac{2}{3} \,\ve}  \stackrel{(\ref{3.2})}{\ge}  p_0(u_0) \,.
\end{equation*}

\n
This completes the proof of (\ref{3.6}). We can thus select $L_0 = c(u)$, such that with the choices (\ref{3.3})-(\ref{3.5}), $u_\infty < u$, and the assumptions of Proposition \ref{prop2.5} are satisfied, hence:
\begin{equation}\label{3.7}
p_n(u) \le p_n(u_n) \le e^{-(K_0 - \log 2) 2^n} \;\mbox{for all $n \ge 0$} \,.
\end{equation}

\n
Together with (\ref{2.11}) and (\ref{2.13}) this implies that for all $n \ge 0$ and $m \in I_n$,
\begin{equation}\label{3.8}
 \IP[C_m \stackrel{\cV^u}{\longleftrightarrow} \wt{S}_m]   \stackrel{(\ref{2.13}}{\le} |\Lambda_m| p_n(u) 
 \underset{(\ref{3.7})}{\stackrel{(\ref{2.11})}{\le}} (c_0 \ell_0^{2(d-1)})^{2^n-1} (2 c_0 \,\ell_0^{2(d-1)})^{-2^n} \le 2^{-2^n}\,.
\end{equation}
Observe that setting $\rho = \frac{log 2}{\log \ell_0}$, we have $2^n = (\frac{L_n}{L_0})^\rho$. Hence for $x$ in $\IZ^d$ outside $B(0,2L_0)$, we pick $L_n$ such that $2L_n < |x|_\infty \le 2L_{n+1}$, and find with (\ref{3.8}) that
\begin{equation}\label{3.9}
\IP [0 \stackrel{\cV^u}{\longleftrightarrow} x] \le c(u) e^{-c^\prime(u) |x|^\rho_\infty}, \;\mbox{for $x$ with $|x|_\infty > 2L_0$}\,.
\end{equation}

\n
Adjusting constants we can ensure that this holds when $|x|_\infty \le 2 L_0$ as well. The claim (\ref{0.3}) now readily follows, and Theorem \ref{theo0.1} is proven. \hfill $\square$

\begin{remark}\label{rem3.1} ~\rm

\medskip\n
1) The results in the present note do not settle the important question of knowing whether $u_*$ and $u_{**}$ coincide. Let us mention that in case they differ our results point to a marked transition in the decay properties of $\IP[0 \stackrel{\cV^u}{\longleftrightarrow} S(0,L)]$, for large $L$ and $u > u_*$, depending on whether $u_* < u < u_{**}$ or $u > u_{**}$. Indeed from the absence of an infinite cluster in $\cV^u$ one knows that
\begin{equation*}
\lim\limits_{L \r \infty} \;\IP[0 \stackrel{\cV^u}{\longleftrightarrow} S(0,L)] = 0, \;\mbox{for $u > u_*$}\,.
\end{equation*}
Then it follows from (\ref{0.2}) that
\begin{equation}\label{3.10}
\mbox{when $u_* < u < u_{**}$, $\underset{L \r \infty}{\overline{\lim}} \;L^{(d-1) + \eta} \,\IP[0 \stackrel{\cV^u}{\longleftrightarrow} S(0,L)] = \infty$, for all $\eta > 0$},
\end{equation}
whereas in view of (\ref{0.3}), one finds that
\begin{equation}\label{3.11}
\mbox{for $u >  u_{**}, \; \IP[0 \stackrel{\cV^u}{\longleftrightarrow} S(0,L)]$ has a stretched exponential decay in $L$}.
\end{equation}

\n
2) One can replace $\cI^u$ and $\cV^u$ with $\cI^{u,R}$, the closed $R$-neighborhood for the $\ell^\infty$-distance of $\cI^u$, and $\cV^{u,R}$, its complement in $\IZ^d$. The proof of Theorem \ref{theo0.1} with minor modifications now shows that for any $u>0$, there exist a positive integer $R$, and $c,c^\prime, \delta >0$, all depending on $d$ and $u$, for which one has (the notation is as in (\ref{0.3})):
\begin{equation}\label{3.12} 
\IP[0 \stackrel{\cV^{u,R}}{\longleftrightarrow} x] \le \ c \, \exp\{-c^\prime |x|^\delta\}, \;\mbox{for all $x \in \IZ^d$}\,.
\end{equation}

\n
The proof of (\ref{3.12}) involves the following changes. In (\ref{2.1}) one now assumes that $L_0 >R$. One defines $A^{u,R}_m$ as $A^u_m$ in (\ref{2.12}), with $\cV^u$ replaced by  $\cV^{u,R}$, and $\wt{S}_m$ replaced by $\wt{S}^R_m$, the set of points of $\wt{C}_m$ at $\ell^\infty$-distance $R$ from $\wt{S}_m$. Introducing $p^R_n(u)$ as $p_n(u)$ in (\ref{2.13}), with $A^u_\cT$ replaced by $A^{u,R}_\cT$, (defined similarly as $A^u_\cT$, with $A^{u,R}_{m^\prime}$ in place of $A^u_{m^\prime}$, $m^\prime \in \cT^0$), Lemma \ref{lem2.1} and Proposition \ref{prop2.2} hold with $A^{u,R}_m$ and $p^R_n(u)$ in place of $A^u_m$ and $p_n(u)$. For the proof of Proposition \ref{prop2.2} one only needs to replace in (\ref{2.37}) and in the definition of $\cI^*$ below (\ref{2.37}), $w(\IN)$ as well as range $w_i$, $1\le i \le \ell$, by their closed $R$-neighborhood for the $\ell^\infty$-distance. Then Proposition \ref{prop2.5}  holds for $p^R_n(u)$ with the additional constraint $L_0 >R$, and $p^R_n(u_0)$ in place of $p_n(u_0)$ in (\ref{2.46}).  In order to initiate the induction and check (\ref{2.45}), (\ref{2.46}), one picks $u_0 = \frac{1}{2}u$, $r_0 = 12$, $K_0$ as in (\ref{3.4}), $\ell_0 = L_0$, and $R = [L_0^{\frac{3}{2(d-1)}}] +1$. When $L_0$ is large, one readily sees that (\ref{2.45}) holds. To check (\ref{2.46}) one notes that when the closed $R$-neighborhood for the $\ell^\infty$-distance of the starting point of the trajectories in the support of  $\mu_{C_m,u_0}$ covers $S_m$, then $A^{u,R}_m$ does not hold.  One further notes with a straightforward lower bound on $P_z[\wt{H}_{C_m} = \infty]$, for $z \in S_m$, that in the notation of (\ref{1.2}),  $e_{C_m}(B(z,R)) \ge c R^{d-1} L_0^{-1}$, for  $z \in S_m$ and $m \in I_0$. As a result one finds that 
\begin{equation*}
p^R_0(u_0) \le c L_0^{d-1} \exp\{-u_0 c R^{d-1}L_0^{-1}\}  \le c L_0^{d-1} \exp\{-c u L_0^{\frac{1}{2}}\} .
\end{equation*}
Hence (\ref{2.46}) holds when $L_0 \ge c(u)$. The proof of (\ref{3.12}) then proceeds in the same way as the proof of (\ref{0.3}).

\n
3) It is an interesting problem to determine the exact nature of the decay in $L$ of $\IP[0 \stackrel{\cV^u}{\longleftrightarrow} S(0,L)]$, for $u > u_*$. Let us point out that when $d=3$, this decay cannot be exponential, as in the case of subcritical Bernoulli percolation, cf.~\cite{Grim99} Theorem 5.4, p.~86. Indeed denoting with $S_L$ the discrete segment along the first coordinate positive half-axis joining $0$ with $S(0,L)$ one has for any $u \ge 0$,
\begin{equation*}
\begin{split}
\IP [0  \stackrel{\cV^u}{\longleftrightarrow} S(0,L)] \ge \IP[\cV^u \supseteq S_L] & \stackrel{(\ref{1.7})}{=} \exp\{ - u \,{\rm cap} (S_L)\}
\\[1ex]
&\; \ge\; \,\exp \Big\{- c\,u \;\mbox{\f $\dis\frac{L}{\log L}$} \Big\}, \;\mbox{for $L \ge 1$, when $d=3$}\,,
\end{split}
\end{equation*}

\n
using for instance (\ref{1.4}), (\ref{1.5}) and standard estimates on the Green function to bound ${\rm cap}(S_L)$. When $d \ge 4$, ${\rm cap}(S_L)$ grows linearly with $L$, and the above calculation does not preclude an exponential decay of $\IP [0  \stackrel{\cV^u}{\longleftrightarrow} S(0,L)] $ in $L$ for some (or all) $u > u_*$. \hfill $\square$
\end{remark}

\end{document}